\documentclass{jancl}

\journal{13}{1}{2003}{-}{-}

\usepackage[all]{xy}

\renewcommand{\b}{\beta}

\newcommand{\G}{\Gamma}

\renewcommand{\t}{\tau}

\renewcommand{\l}{\lambda}
\renewcommand{\L}{\Lambda}

\newcommand{\sou}{\overline}

\newcommand{\f}{\rightarrow}
\newcommand{\q}{\forall}

\renewcommand{\v}{\vdash}

\newcommand{\p}{\succ}

\newenvironment{The}{\begin{theorem}}{\end{theorem}}
\newenvironment{Pro}{\begin{proposition}}{\end{proposition}}
\newenvironment{Cor}{\begin{corollary}}{\end{corollary}}
\newenvironment{Lem}{\begin{lemma}}{\end{lemma}}

\title[Complete Types]{Complete Types in an Extension of the \\System ${\cal AF}2$}

\author{Samir Farkh \andauthor Karim Nour}

\address{
LAMA - \'Equipe de Logique\\
Université de Chambéry\\
73376 Le Bourget du Lac\\[3pt]
nour@univ-savoie.fr
}

\abstract{In this paper, we extend the system ${\cal AF}2$ in
order to have the subject reduction for the $\beta
\eta$-reduction. We prove that the types with positive quantifiers
are complete for models that are stable by weak-head expansion.}


\begin{document}

\maketitle

\section{Introduction}

The semantics of realisability of the system ${\cal F}$, proposed
by J.-Y.\ Girard, consists in interpreting the types by
``saturated subsets'' of $\l$-terms. The correction theorem (also
called ``adequacy lemma'') stipulates that: if a $\l$-term is
typable then it belongs to the interpretation of its type. The
adequacy lemma allows to show the strong normalization of the
system ${\cal F}$ when we take an adequate concept of saturation.
The power of this notion of semantics comes from the variety of
possible interpretations of the second order quantifier.  For the
system ${{\cal AF}2}$, J.-L.\ Krivine proposed a more general
semantics by defining the concept of $\l$-models for a
second-order language. His semantics is a modification of the
traditional concept of a second-order model in which the set of
the truth values is not, as usual, $\{0, 1\}$ but an adequate
subset of $\l$-terms (see \cite{KRI:1994} and \cite{RAF:1998}).
The corresponding adequacy lemma allows also to prove the
uniqueness of the representation of the data.

Many researchers were interested in finding a general definition
of a data type. For example, Böhm and Berarducci gave such a
definition, only for term algebras, in the system ${\cal F}$ (see
\cite{BB:1985}) and Krivine generalized their definition to system
${{\cal AF}2}$ (see \cite{KRI:1990b}). We noticed that the class
${\cal A}$ of the types thus built has the following feature: a
normal $\l$-term is typable of a type $D \in { \cal A}$ iff it is
in the interpretation of $D$ for a certain semantics. Then, we
decided to take this result as the definition of the data types
which we called ``complete types'', because the considered
semantics is complete for these types. R.\ Labib-Sami was the
first to build a class of complete types: they are the types with
positive quantifiers (denoted by $\q_2^+$) of the system ${\cal
F}$ compared to a semantics based on the subsets saturated by $\b
\eta$-equivalence (see \cite{LAB:1986}).

\begin{sloppypar}
We generalized in \cite{FN:1998b} Labib-Sami's result, by showing
that the $\q_2^+$ types of system ${{\cal AF}2}$ are complete for
the semantics based on sets saturated by $\b \eta$-equivalence.
It was natural to imagine a refinement of this result, namely
interpretation of the types by sets saturated by weak-head
expansion.  For this, we considered a more restricted class of the
$\q_2^+$ types which includes the data types of J.-L.\ Krivine.
Then, we showed in  \cite{FN:1998b} that these types are preserved
by $\eta$-reduction and are complete for the considered semantics.
\end{sloppypar}

We propose in this paper another solution to this problem. We add
typing rules to the system ${{\cal AF}2}$ in order to have the
conservation of types by $\b \eta$-reduction. The system, which we
propose, is inspired by the works of Mitchell \cite{MIT:1988} and
the second author \cite{NOUR:1996}. We show that, in this new
system, all $\q_2^+$ types are complete for the semantics based on
sets saturated by weak-head expansion.

\section{Notations and definitions}

\begin{notations}
We denote by $\L$ the set of terms of pure $\l$-calculus, also
called $\l$-terms. Let $t,u,u_1,\dots,u_n \in \L$, the application
of $t$ to $u$ is denoted by $(t)u$. In the same way we write
$(t)u_1\ldots u_n$ instead of $(\ldots((t)u_1)\ldots)u_n$. The
$\b$-reduction (resp.\ $\b$-equivalence) is denoted by $t \f_{\b}
u$ (resp.\ $t \simeq\sb{\b} u$). The set of free variables of a
$\l$-term $t$ is denoted by $Fv(t)$. Let us recall that a
$\l$-term $t$ either has a \textbf{weak-head redex} [i.e.\ $t=(\l
xu)vv_1 \ldots v_m$, the weak-head redex being $(\l xu)v$], or is
in \textbf{weak-head normal form} [i.e.\ $t=(x)vv_1\ldots v_m$ or
$t= \l x u$].  The notation $t \p_f t'$ means that $t'$ is
obtained from $t$ by some \textbf{weak-head reductions}.
\end{notations}

\subsection{The ${{\cal AF}2}$ type system}

The types will be formulas of second-order predicate logic over a
given language. The logical symbols are $\perp$ (for absurd), $\f$
and $\q$ (and no other ones). There are individual variables:
$x,y,\ldots$ (also called first-order variables) and $n$-ary
predicate variables ($n=0,1,\ldots$): $X,Y,\ldots$ (also called
second-order variables). The terms and formulas are built in the
usual way.

If $X$ is a unary predicate variable, $t$ and $t'$ two terms, then
the formula $\q X[Xt \f Xt']$ is denoted by $t=t'$, and is said to
be an \textbf{equation}. A \textbf{particular case} of $t=t'$ is a
formula of the form
$t[u_1/x_1,\ldots,u_n/x_n]=t'[u_1/x_1,\ldots,u_n/x_n]$ or
$t'[u_1/x_1,\ldots,u_n/x_n]=t[u_1/x_1,\ldots,u_n/x_n]$,
$u_1,\ldots,u_n$ being terms of the language. Then, we denote by
\textbf{E}  a system of function equations. A context $\G$ is a
set of the form $x_1:A_1,\ldots,x_n:A_n$ where $x_1,\ldots,x_n$
are distinct variables and $A_1,\ldots,A_n$ are formulas. We are
going to describe a system of typed $\l$-calculus called
second-order functional arithmetic (abbreviated in ${\cal AF}2$
for \textbf{Arithmétique Fonctionnelle du second ordre}). The
typing rules are the following:
\begin{janclenum}{(\arabic}{)}
\item $\G,x:A \v_{{\cal AF}2} x:A$.

\item If $\G,x:B \v_{{\cal AF}2} t:C$, then $\G \v_{{\cal AF}2}\l
xt:B \f C$.

\item If $\G \v_{{\cal AF}2}u:B \f C$, and $\G \v_{AF2}v:B$, then
$\G\v_{{\cal AF}2} (u)v:C$.

\item If $\G \v_{{\cal AF}2}t:A$, and $x$ does not appear in $\G$,
then $\G \v_{{\cal AF}2} t:\q xA$.

\item If $\G \v_{{\cal AF}2}t:\q xA$, then, for every term $u$,
$\G \v_{{\cal AF}2}t:A[u/x]$.

\item If $\G \v_{{\cal AF}2}t:A$, and $X$ does not appear in $\G$,
then $\G \v_{{\cal AF}2}t:\q XA$.

\item If $\G \v_{{\cal AF}2}t:\q XA$, then, for every formula $G$,
\item[] $\G \v_{{\cal
AF}2}t:A[G/X(x_1,\ldots,x_n)]$.\footnote{$A[G/X(x_1,\ldots,x_n)]$
is obtained by replacing in $A$ each atomic formula
$X(t_1,\ldots,t_n)$ by $G[t_1/x_1,\ldots,t_n/x_n]$. To simplify,
we write $A[G/X]$ instead of $A[G/X(x_1,\ldots,x_n)]$.}

\item If $\G \v_{{\cal AF}2}t:A[u/x]$, then $\G \v_{{\cal
AF}2}t:A[v/x]$, $u=v$ being a particular case of an equation of
\textbf{E}.
\end{janclenum}

Whenever we obtain the typing $\G \v_{{\cal AF}2}t:A$ by means of
these rules, we say that ``the $\l$-term $t$ is of type $A$ in the
context $\G$, with respect to the equations of \textbf{E}''.

\begin{The}\label{theorem:un}
\begin{enumerate}
\item If $\G \v_{{\cal AF}2}t:A$, and $t \f_{\b} t'$, then $\G
\v_{{\cal AF}2}t':A$.

\item If $\G \v_{{\cal AF}2}t:A$, then $t$ is strongly
normalizable.
\end{enumerate}
\end{The}

\subsection{The semantics of ${{\cal AF}2}$}

 If $G,G' \in P(\L)$, we define an element of $P(\L)$ by: $G\f G'
= \{ u\in \L$ / $(u)t \in G'$ for every $t\in G\}$.
 Let ${\cal R}_f$ the set of subsets of $\L$ stable by
weak-head reduction (i.e.\ $\Xi \in {\cal R}_f$ iff for every $v
\in \Xi$, if $u \p_f v$, then $u \in \Xi$). A subset $R$ of ${\cal
R}_f$ is said \textbf{adequate} iff $R$ is closed by $\f$ and
$\cap$.

Let $L$ be a second-order language. A \textbf{$\L_f$-model} is
defined by:
\begin{itemize}
\item a non empty set $|M|$ called \textbf{domain} of $M$,

\item an adequate set $R$ of ${\cal R}_f$,

\item for every an $n$-ary function symbol of $L$,
 a function $f_M :  |M|^n \f |M|$,

\item for every $n$-ary predicate symbol $P$ of $L$, a function
$P_M :  | M |^n \f R$.
\end{itemize}

Let $M$ be a $\L$-model of $L$.
\begin{itemize}
\item An \textbf{interpretation} $I$ is a function from the set of
first (resp.\ the set of $n$-ary second) order variables to $|M|$
(resp.\ to $R^{ | M |^n}$).

\item Let $I$ be an interpretation, $x$ (resp.\ $X$) a first
(resp.\ an $n$-ary second) order variable, and  $a$ (resp.\
$\Phi$) an element of $|M|$ (resp.\ of $R^{ |M|^n}$).  We define
an interpretation $J = I [x \leftarrow a]$ (resp.\ $J = I [X
\leftarrow \Phi]$) by taking $J (x) = a$ (resp.\ $J (X) = \Phi$)
and $J (\xi) = I (\xi)$ (resp.\ $J (\xi') = I (\xi')$) for every
variables $\xi \not = x$ (resp.\ $\xi' \not = X$).
\end{itemize}

Let $I$ be an interpretation. To every term $t$ of $L$, we define,
by induction, its \textbf{value} $t_{M,I} \in |M|$:
\begin{itemize}
\item if $t=x$, then $t_{M,I} = I (x)$,

\item if $t = f(t^1,\ldots, t^n)$, then $t_{M,I} = f_M
(t_{M,I}^1,\ldots, t_{M,I}^n)$.
\end {itemize}

Let $A$ be a formula of $L$. The \textbf{value} of $A$ in a model
$M$ and an interpretation $I$ (denoted by \textbf{$|A|_{M,I}$}) is
an element of $R$ defined by induction:
\begin{itemize}
\item if $A=P(t^1,\ldots, t^n)$, where $P$ is a $n$-ary predicate
symbol
  (resp.\ second-order variable) and  $t^1,\ldots, t^n$ are terms of $L$,
  then $|A|_{M,I} = P_M (t_{M,I}^1,\ldots, t_{M,I}^n)$ (resp.\ $|A|_{M,I}
  = I (X)(t_{M,I}^1,\ldots,t_{M,I}^n)$).

\item if $A = B\f C$, then $|A|_{M,I} = |B|_{M,I} \f |C|_{M,I}$,

\item if $A = \q xB$ where $x$ is a first-order variable, then
$|A|_{M,I} = \bigcap \{|B [x]|_{M,I [x \leftarrow a]}$; $a\in | M
|\}$,

\item if $A = \q XB$ where $X$ is a $n$-ary second-order variable,
then $|A|_{M,I} = \\ \bigcap \{ |B[X]|_{M,I [X \leftarrow \Phi]}$;
$\Phi \in R^{|M|^n}\}$.
\end {itemize}

It is clear that: if $A$ is a closed type, then $|A|_{M,I}$ does not
depend on the interpretation $I$ and we write $|A|_M$.

Let $M$ be a $\L$-model of $L$.
\begin{itemize}
\item We say that \textbf{$M$ satisfies} the equation $u = v$, if
for every interpretation $I$, $u_{M,I} = v_{M,I}$. If \textbf{$E$}
is a set of equations of $L$, we say that $M$ \textbf{satisfies}
$E$, or $M$ is a \textbf{model} for $E$, iff $M$ satisfies all the
equations of  $E$.

\item If $A$ is a closed formula, we denote by $|A|_f = \bigcap \{
|A]|_{M}$; $M$ is a $\L_f$-model which satisfies $E\}$.
\end{itemize}


The following theorem is  known under the name ``adequation lemma'' or
``the correction theorem'':
\begin{The}\label{theorem:deux}
Let $t$ be a $\l$-term and $A$ a closed type of system ${{\cal
AF}2}$. If $\v_{{\cal AF}2} t : A$, then $t \in |A|_f$.
\end{The}

\section{The system ${{\cal AF}2_ {\subseteq}}$}

The typing system ${\cal AF}2$ does not conserve the types by
$\eta$-reduction. Indeed, $\v_{{\cal AF}2} \l x \l y (x)y : \q X
(X \f (X \f X)) \f (\q XX \f \q X(X \f X))$ but $\l x \l y (x)y
\f_{\eta} \l xx$ and $\not \v_{{\cal AF}2} \l xx : \q X (X \f (X
\f X)) \f (\q XX \f \q X(X \f X))$. We will define an extension of
the system ${\cal AF}2$ which, while keeping the properties of the
system ${\cal AF}2$, conserves the types by $\eta$-reduction.

\begin{definition}
Let $E$ be an equation system of second-order language $L$.  We
define on the formulas of ${{\cal AF}2}$ a binary relation
$\subseteq$ by: $A \subseteq B$ iff it is obtained by the
following proof rules:
 \begin{description}
\item[$(ax')$]  $A \subseteq A$
\item[$(dist)$] $\q \xi (C \f D) \subseteq \q \xi C \f \q \xi D$
\item[$(\f)$] If $C' \subseteq C$ and $ D \subseteq D'$, then $C \f D \subseteq C' \f D'$
\item[$(\q_e)$] If $A \subseteq \q \xi C$, then  $A \subseteq C [F/\xi]$
\item[$(\q_i)$] If $A \subseteq D$ and $\xi$ is not free in  $A$, then  $A \subseteq \q \xi D$
\item[$(tr)$] If $A \subseteq D$ and $D \subseteq B$, then  $A \subseteq B$
\item[$(e)$] If $A \subseteq D [u/y]$ and $u = v$ is a particular case
  of an equation of $E$, then  $A \subseteq D [v/y]$
\end{description}
\end{definition}

\begin{definition}
The system ${{\cal AF}2_{\subseteq}}$ is the system ${{\cal AF}2}$
 where we add the following rule:
\begin{gather}
\text{If }\G \v_{{\cal AF}2_ {\subseteq}} t : A\text{ and  }A
\subseteq B\text{, then }\G \v_{{\cal AF}2_{\subseteq}} t : B
\tag{$\subseteq$}\label{rule:subseteq}
\end{gather}
\end{definition}
It is clear that the rules $(5)$, $(7)$ and $(8)$ are
particular cases of the rule ($\subseteq$).

\subsection{Syntactical properties of the system}

\begin{notations} Let $\mbox{\boldmath$\xi$} = \xi_1,\ldots, \xi_n$ be
a sequence of variables.  We denote the formula $\q \xi_1\ldots \q
\xi_nF$ by $\q \mbox{\boldmath$\xi$} F$. We write
``\mbox{\boldmath$\xi$} is not free in $A$'' if for every $1 \leq
i \leq n$, $\xi_i$ is not free in $A$.  Let $A$ be a formula,
\textbf{F} a sequence of formulas $F_1,\ldots, F_n$,
\mbox{\boldmath$u$} a sequence of terms $u_1,\ldots, u_n$ and
\mbox{\boldmath$x$} (resp.\ \mbox{\boldmath$X$}) a sequence of
first (resp second) ordre variables $x_1,\ldots, x_n$ (resp.\
$X_1,\ldots, X_n$). We denote by $A
[\mbox{\boldmath$u$}/\mbox{\boldmath$x$}]$ the formula $A
[u_1/x_1,\ldots, u_n/x_n]$ and by
$A[\textbf{F}/\mbox{\boldmath$X$}]$
 the formula $A [F_1/X_1,\ldots, F_n/X_n]$.
 \end{notations}

\begin{Lem}\label{lemme:un}
In the typing, we may replace the succession of $n$ times
($\subseteq$) and $m$ times $(4)$ and $(6)$, by the
succession of $m$ times $(4)$ and $(6)$, and $n$ times
($\subseteq$).
\end{Lem}

\begin{proof*}
By induction on $n$ and $m$.
\end{proof*}

We deduce the following corollary:
\begin{Cor}\label{corollaire:un}
If  $\G \v_{{\cal AF}2_{ \subseteq}} t : B$ is derived from $\G
\v_{{\cal AF}2_{\subseteq}} t : A$, then we may assume that we
begin by the applications of $(4)$, $(6)$ and next
($\subseteq$) (i.e.\ there is  \mbox{\boldmath$\xi$} not free in
$\G$ such that $\q \mbox{\boldmath$\xi$} A  \subseteq B$).
\end{Cor}

Then we have the following  characterization:
\begin{The}\label{theorem:trois}
\begin{janclenum}{(\roman}{)}%
\item If $\G \v_{{\cal AF}2_{\subseteq}} x : A$, then there
  is a type $B$ such that $x : B \in \G$ and $\q \mbox{\boldmath$\xi$}
   B \subseteq A$, where \mbox{\boldmath$\xi$} is not free in $\G$.
\item If $\G \v_{{\cal AF}2_{\subseteq}} \l xu : A$, then there
are two types $B$ and $C$ such that $\G, x : B \v_{{\cal
AF}2_{\subseteq}} u : C$ and $\q \mbox{\boldmath$\xi$} (B \f C)
\subseteq A$, where \mbox{\boldmath$\xi$} is not free in $\G$.

\item If $\G \v_{{\cal AF}2_{\subseteq}} (u)v : A$, then there are
two types $B$ and $C$ such that $\G \v_{{\cal AF}2_{\subseteq}} u
: B \f C$, $\G \v_{{\cal AF}2_{\subseteq}} v : B$ and $\q
\mbox{\boldmath$\xi$} C \subseteq A$, where \mbox{\boldmath$\xi$}
is not free in $\G$.
\end{janclenum}
\end{The}

 We will define a typing system equivalent to system ${{\cal AF}2_{\subseteq}}$.
\begin{definition} The system ${{\cal AF}2S}$ is defined only
 by the three following rules:
\begin{description}
\item[(S1)] If $x : B \in \G$ and $\q \mbox{\boldmath$\xi$} B
 \subseteq A$, then $\G \v_{{\cal AF}2S} x : A$
\item[(S2)] If $\G, x : B \v_{{\cal AF}2S} u : C$ and $ \q \mbox{\boldmath$\xi$} (B \f C)
 \subseteq A$, then $\G\v_{{\cal AF}2S} \l xu : A$
\item[(S3)] If $\G\v_{{\cal AF}2S} u : B \f C$, $\G\v_{{\cal AF}2S}
 v : B$ and $\q \mbox{\boldmath$\xi$} C \subseteq A$, then  $\G\v_{{\cal AF}2S} (u)v : A$
\end{description}
where \mbox{\boldmath$\xi$} is not free in $\G$.
\end{definition}

We have the following result:
\begin{The}\label{theorem:quatre}
$\G \v_{{\cal AF}2_{\subseteq}} t : A$ iff $\G \v_{{\cal AF}2S} t : A$.
\end{The}

\begin{proof*}
We use Theorem \ref{theorem:trois}.
\end{proof*}

In the rest of the paper we often consider the system ${{\cal AF}2S}$.\\

The following corollary will often be used:
\begin{Cor}\label{corollaire:deux}
If $\G, x : A \v_{{\cal AF}2S} (x)u_1\ldots u_n : B$, then  \\
$n = 0$, $\q \mbox{\boldmath$\xi_0$} A  \subseteq B$ and
$\mbox{\boldmath$\xi_0$}$ does not appear in $\G$ and $A$, or\\
$n \neq 0$,   $\q \mbox{\boldmath$\xi_0$} A  \subseteq C_1 \f
B_1$, $\q \mbox{\boldmath$\xi_i$} B_i  \subseteq  C_{i + 1} \f
B_{i + 1}$ ($1 \leq i \leq n-1)$, and $\q \mbox{\boldmath$\xi_n$}
B_n  \subseteq B$ where \mbox{\boldmath$\xi_i$} ($0 \leq i \leq
n)$ are not free in $\G$ and $A$, and $\G, x : A \v_{{\cal AF}2S}
u_i : C_i$ ($1 \leq i \leq n)$.
\end{Cor}

\begin{proof*}
By induction on $n$ and using Theorem \ref{theorem:trois}.
\end{proof*}

\subsection{Conservation of type by $\b$-reduction}

\begin{Lem}\label{lemme:deux}
If $A \subseteq B$, then, for any sequence of terms $\textbf{u}$
(resp.\ of formulas $\textbf{F}$), $A[\textbf{u}/\textbf{x}]
\subseteq B [\textbf{u}/\textbf{x}]$ (resp.\
$A[\textbf{F}/\textbf{X}] \subseteq B [\textbf{ F}/\textbf{X}]$),
and we use the same proof rules.
\end{Lem}

\begin{proof*}
By induction on the derivation  $A \subseteq B$.
\end{proof*}

\begin{Lem}\label{lemme:trois}
If $\G \v_{{{\cal AF}2S}} t : A$, then, for all sequences of terms
$\textbf{u}$ (resp.\ of formulas $\textbf{F}$), $\G
[\textbf{u}/\textbf{x}] \v_{{{\cal AF}2S}} t : A [\textbf{
u}/\textbf{x}]$ (resp.\ $\G [\textbf{F}/\textbf{Y}]\v_{{{\cal
AF}2S}} t :
 A[\textbf{F}/\textbf{Y}]$) and we use the same typing rules.
\end{Lem}

\begin{proof*}
By induction on the derivation  $\G \v_{{{\cal AF}2S}} t : A$. We
look at the last rule used and we use Lemma \ref{lemme:deux}.
\end{proof*}

\begin{Lem}\label{lemme:quatre}
If $x_1 : A_1,\ldots, x_n : A_n \v_{{{\cal AF}2S}} t : A$, $B_i
\subseteq A_i$ ($1 \leq i \leq n$) et $A \subseteq B$, then $x_1 :
B_1,\ldots, x_n : B_n \v_{{{\cal AF}2S}} t : B$.
\end{Lem}

\begin{proof*}
By induction on the $\l$-term $t$.
\end{proof*}

\begin{Lem}\label{lemme:cinq}
If $\G, x : B \v_{{{\cal AF}2S}} u : A$ et $\G \v_{{{\cal AF}2S}}
v : B$, then $\G \v_{{{\cal AF}2S}} u [v/x] : A$.
\end{Lem}

\begin{proof*}
By induction on the derivation $\G, x : B \v_{{{\cal AF}2S}} u :
A$.
\end{proof*}

\begin{Lem}\label{lemme:six}
If $\G, x : C \v_{{{\cal AF}2S}} u : D$ and there is a
\mbox{\boldmath$\xi$} which does not appear in $\G$ and
\mbox{\boldmath$\xi'$} such that $\q \mbox{\boldmath$\xi$} (C\f D)
\subseteq \q \mbox{\boldmath$\xi'$} (A\f B)$, then $\G, x : A
\v_{{{\cal AF}2S}} u : B$.
\end{Lem}

\begin{proof*}
By induction on the derivation $\q \mbox{\boldmath$\xi$} (C\f D)
\subseteq \q \mbox{\boldmath$\xi'$} (A\f B)$. We look at the last
rule used. We consider only three cases.
\begin{description}
\item[$(\f)$] We have $A \subseteq C$ and $D \subseteq B$,
 then, by Lemma \ref{lemme:quatre}, we deduce the result.

\item[$(e)$] We have $\q \mbox{\boldmath$\xi$} (C\f D) \subseteq E
[u/y] = \q \mbox{\boldmath$\xi'$} (F [u/y] \f M [u/y])$. Then $A =
F [v/y]$ and $B = M [v/y]$ where $u = v$ is a particular case of
an equation of $E$. By induction hypothesis, we obtain $\G, x : F
[u/y] \v_{{{\cal AF}2S}} u : M [u/y]$. But $F [v/y] \subseteq F
[u/y]$ and $M [u/y] \subseteq M [v/y]$, then, by Lemma
\ref{lemme:quatre}, $\G, x : A \v_{{{\cal AF}2S}} u : B$.

\item[$(\q_e)$] We have $\q \mbox{\boldmath$\xi$} (C\f D)
 \subseteq \q s \q \mbox{\boldmath$\xi'$} (E\f F)$ and $A = E [ G/s]$,
 $B = F [G/s]$. By induction hypothesis, we obtain $\G, x : E
 \v_{{{\cal AF}2S}} u : F$. We may assume that $s$ is not free in
 $\G$, then, by Lemma \ref{lemme:trois}, $\G, x : E [ G/s] \v_{{{\cal AF}2S}} u : F
 [G/s]$, i.e $\G, x : A \v_{{{\cal AF}2S}} u : B$.
\end{description}
\end{proof*}

\begin{Lem}\label{lemme:sept}
If $\G \v_{{{\cal AF}2S}} \l xu : A\f B$, then $\G, x : A
\v_{{{\cal AF}2S}} u : B$.
\end{Lem}

\begin{proof*}
We have $\G \v_{{{\cal AF}2S}} \l xu : A\f B$, then $\G, x : C
\v_{{{\cal AF}2S}} u : D$ and $\q \mbox{\boldmath$\xi$} (C\f D)
\subseteq (A\f B)$ where \mbox{\boldmath$\xi$} is not free in
$\G$. Therefore, by Lemma \ref{lemme:six}, $\G, x : A \v_{{{\cal
AF}2S}} u : B$.
 \end{proof*}

\begin{The}\label{theorem:cinq}
If $\G \v_{{{\cal AF}2S}} t : A$ and $t \f_{\b} t'$, then $\G
\v_{{{\cal AF}2S}} t' : A$.
\end{The}

\begin{proof*}
It suffices to do the proof for one step of reduction. We proceed
by induction on $t$ et we use Lemmas \ref{lemme:cinq} and
\ref{lemme:sept}.
\end{proof*}

\subsection{Conservation of type by $\eta$-reduction}

\begin{The}\label{theorem:six}
If $\G \v_{{\cal AF}2S} t : A$ and $t \f_{\eta} t'$, then $\G
\v_{{\cal AF}2S} t' : A$.
\end{The}

\begin{proof*}
It suffices to do the proof for one step of $\eta$-reduction
denoted $\eta_0$. We do the proof by induction on $t$. The only
difficult case is $t = \l xu$, then two cases can arise:
\begin{enumerate}
\item $t' = \l xu'$ where $u \f_{\eta_0} u'$: We have $\G \v_{{{\cal
 AF}2S}} \l xu : A$, then $\G, x : B \v_{{{\cal AF}2S}} u : C$ and $\q
 \mbox{\boldmath$\xi$}(B\f C ) \subseteq A$ where
 $\mbox{\boldmath$\xi$}$ is not free in $\G$. By induction hypothesis,
 we have $\G, x : B \v_{{{\cal AF}2S}} u' : C$, and, by the rule
 $(S2)$, $\G \v_{{{\cal AF}2S}} \l xu' : A$, i.e $\G \v_{{{\cal
 AF}2S}} t' : A$.

\item $u = (t')x$ where $x$ is not free in $t'$: We have $\G, x :
B \v_{{{\cal AF}2S}} (t')x : C$, and $\q \mbox{\boldmath$\xi$}(B\f
C ) \subseteq A$ where $\mbox{\boldmath$\xi$}$ is not free in
$\G$. Then $\G, x : B \v_{{{\cal AF}2S}} t' : E\f F$, $\G, x : B
\v_{{{\cal AF}2S}} x : E$ and $\q \mbox{\boldmath$\xi'$} F
\subseteq  C$ where $\mbox{\boldmath$\xi'$}$ is not free in $\G$
and $B$. By  Corollary \ref{corollaire:deux}, we obtain $\q
\mbox{\boldmath$\xi''$} B \subseteq E$ where
$\mbox{\boldmath$\xi''$}$ is not free in $\G$ and $B$. We have $B
\subseteq B$, then $B \subseteq \q \mbox{\boldmath$\xi''$}B
\subseteq E$, and $B \subseteq \q \mbox{\boldmath$\xi'$}E$.  Using
the rules $(dist)$ and $(\f)$, we deduce $\q
\mbox{\boldmath$\xi'$} (E\f F) \subseteq \q \mbox{\boldmath$\xi'$}
E \f \q \mbox{\boldmath$\xi'$} F \subseteq B\f C$ and $\q
\mbox{\boldmath$\xi$} \q \mbox{\boldmath$\xi'$} (E\f F) \subseteq
\q \mbox{\boldmath$\xi$} (B\f C)$. Finally , we have $\G
\v_{{{\cal AF}2S}} t' : E\f F$, then $\G \v_{{{\cal AF}2S}} t' :
\q \mbox{\boldmath$\xi$} \q \mbox{\boldmath$\xi'$} (E\f F)$, and,
by the rule $(tr)$, we obtain $\G \v_{{{\cal AF}2S}} t' : A$.
\end{enumerate}
\end{proof*}

We will see that the system ${{\cal AF}2S}$ is exactly  ${{\cal
 AF}2}$ in which one adds the conservation of the type by
 $\eta$-reduction as a typing rule.

\begin{definition}
The typing system ${{\cal AF}2{\eta}}$ is the system ${{\cal
AF}2}$, in which we add the following typing rule:
\begin{gather} \text{If }\G \v_{{{\cal AF}2{\eta}}} t : A\text{
and }t \f_{\eta} t'\text{, then }\G \v_{{{\cal AF}2{\eta}}} t' : A
\tag{$\eta$}\label{rule:eta}
\end{gather}
\end{definition}
The typing rule (\ref{rule:subseteq}) is derivable in the system
${{\cal AF}2{\eta}}$.
\begin{The}\label{theorem:sept}
If $\G \v_{{{\cal AF}2{\eta}}} t : A$ and $A \subseteq B$, then
$\G
 \v_{{{\cal AF}2{\eta}}} t : B$.
\end{The}
\begin{proof*}
By induction on the proof of $A \subseteq B$. We consider the last
rule used. The only difficult case is $(dist)$. We have $A =
 \q \xi (C\f D)$ and $B = \q \xi C \f \q \xi D$. If $\G, x : \q
 \mbox{\boldmath$\xi$} C \v_{{{\cal AF}2}} t : \q
 \mbox{\boldmath$\xi$} (C \f D)$, then $\G, x : \q
 \mbox{\boldmath$\xi$} C \v_{{{\cal AF}2}} t : C \f D$ and $\G,
 x : \q \mbox{\boldmath$\xi$} C \v_{{{\cal AF}2}} (t)x : D$.  Since
 \mbox{\boldmath$\xi$} is not free in $\G$, we obtain $\G, x : \q
 \mbox{\boldmath$\xi$} C \v_{{{\cal AF}2}} (t)x : \q
 \mbox{\boldmath$\xi$} D$ and $\G \v_{{{\cal AF}2}} \l x
 (t)x : \q \mbox{\boldmath$\xi$} C \f \q \mbox{\boldmath$\xi$}D$.
 Since $\l x (t)x \f_{\eta} t$, we deduce $\G \v_{{{\cal AF}2}} t :
 B$.
\end{proof*}

We can then deduce the following result:
\begin{The}\label{theorem:huit}
$\G \v_{{{\cal AF}2S}} t : A$ iff  $\G \v_{{{\cal AF}2{\eta}}} t :
A$.
\end{The}

\begin{proof*}
By Theorems \ref{theorem:six} et \ref{theorem:sept}.
\end{proof*}

We can also state the following proposition:
\begin{Pro}\label{proposition:un}
If $\G \v_{{{\cal AF}2{\eta}}} t : A$, then there is a $\l$-term
$u$ such that $u \f_{\eta} t$ and $\G \v_{{{\cal AF}2}} u : A$.
\end{Pro}

\begin{proof*}
By induction on the typing $\G \v_{{{\cal AF}2{\eta}}} t : A$.
\end{proof*}

\subsection{The strong normalization}

\begin{notation}
We write $u\f_{\b^+}v$ if $v$ is obtained from $u$ by at least one
step of $\b$-reduction denoted $\b_0$.
\end{notation}

\begin{Lem}\label{lemme:huit}
Let $u,t,v$ be $\l$-terms such that $u \f_{\eta} t$ and
$t\f_{\b_0}v$. Then there is a $\l$-term $w$ such that $u
\f_{\b^+} w$ and $w \f_{\eta} v$.
\end{Lem}

\begin{proof*}
See \cite{BAR:1984}.
\end{proof*}

\begin{Lem}\label{lemme:neuf}
Let $u, t$ be $\l$-terms. If $u$ is strongly normalizable, and $u
\f_{\eta} t$, then $t$ is also strongly normalizable.
\end{Lem}

\begin{proof*}
If $t$ is not strongly normalizable, then there is an infinite
sequense of $\b_0$-reductions starting with $t$. Since $u
\f_{\eta} t$, then, by Lemma \ref{lemme:huit}, we construct an
infinite sequence of $\b_0$-reductions starting with $u$.
\end{proof*}

\begin{The}\label{theorem:neuf}
If $\G \v_{{\cal AF}2S} t : A$, then $t$ is  strongly
normalizable.
\end{The}

\begin{proof*}
By Proposition \ref{proposition:un}, Theorem \ref{theorem:un} and
Lemma \ref{lemme:neuf}.
\end{proof*}

 \section{The complete types}

\begin{definition}
We say that a closed type $A$ is \textbf{complete} in ${{\cal
AF}2S}$ iff $|A|_{f} = \{ t \in \L$ / $t\f_{\b} t'$ and $\v_{{\cal
AF}2S} t' : A \}$.
\end{definition}

We will give a class of complete types.  We start by extending the
correction theorem to system ${{\cal AF}2S}$.

\begin{Lem}\label{lemme:dix}
Let $M$ be a $\L_f$-model of $E$ and $I$ an interpretation of $E$.
If $A \subseteq B$, then $|A|_{M, I} \subseteq |B|_{M, I}$.
\end{Lem}

\begin{proof*}
By induction on the derivation $A \subseteq B$.
\end{proof*}

\begin{The}[The generalized correction]\label{theorem:dix}
Let $M$ be a $\L_f$-model of $E$ and $I$ an interpretation. If $\G
= x_1 : B_1,\ldots,x_n : B_n \v_{{\cal AF}2S}t' : A$, $t
\simeq\sb{\b} t'$, and $u_i \in |B_i|_{M,I}$ ($1 \leq i \leq n$),
then $t [u_1/x_1,\ldots,u_n/x_n]_{M,I} \in |A|_{M,I}$.
\end{The}

\begin{proof*}
We may assume that $t'$ is normal. The proof is done by induction
on the typing of $t'$. We look at the last rule used.
\begin{description}
\item[$(S1)$] Then $t' = x_i$ ($1\leq i \leq n$) and $\q
\mbox{\boldmath$\xi$} B_i \subseteq A$ where
$\mbox{\boldmath$\xi$}$ is not free in $B_i$ ($1 \leq i \leq n)$.
Since $t \simeq\sb{\b} x_i$, then $t \p_f x_i$ and $t
[u_1/x_1,\ldots,u_n/x_n] \p_f u_i$. But $u_i \in |B_i|_{M,I}$,
then $t [u_1/x_1,\ldots,u_n/x_n] \in |B_i|_{M,I}$. Since
$\mbox{\boldmath$\xi$}$ is not free in $B_i$, we deduce $t
[u_1/x_1,\ldots,u_n/x_n] \in \q \mbox{\boldmath$\xi$} B_i$ and, by
Lemma \ref{lemme:dix}, we obtain $t [u_1/x_1,\ldots,u_n/x_n] \in
|A|_{M,I}$.

\item[$(S2)$] Then $t' = \l xu'$, $\G, x : B \v_{{\cal AF}2S} u' :
C$ and $\q \mbox{\boldmath$\xi$}(B\f C) \subseteq A$ where
$\mbox{\boldmath$\xi$}$ is not free in $B_i$ ($1 \leq i \leq n)$.
Since $t \simeq\sb{\b} \l xu'$, then $t \p_f \l xu$ where $u
\simeq\sb{\b} u'$ and $t [u_1/x_1,\ldots,u_n/x_n] \p_f \l xu
[u_1/x_1,\ldots,u_n/x_n]$. Therefore, by induction hypothesis,
\\$u [u_1/x_1,\ldots,u_n/x_n,v/x] \in |C|_{M,I}$ for all $v \in
|B|_{M,I}$.  We have \\$(\l xu[u_1/x_1,\ldots,u_n/x_n])v \p_f u
[u_1/x_1,\ldots,u_n/x_n,v/x]$, then \\ $\l xu
[u_1/x_1,\ldots,u_n/x_n] \in |B \f C|_{M,I}$, and $t
[u_1/x_1,\ldots,u_n/x_n] \in |\q \mbox{\boldmath$\xi$}(B \f
C)|_{M,I}$. By Lemma \ref{lemme:dix}, we deduce $t
[u_1/x_1,\ldots,u_n/x_n] \in |A|_{M,I}$.

\item[$(S3)$] Then $t' = (u')v'$, $\G \v_{{\cal AF}2S} u' : B\f
C$, $\G \v_{{\cal AF}2S} v' : B$ and $\q \mbox{\boldmath$\xi$} C
\subseteq A$ where $\mbox{\boldmath$\xi$}$ is not free in $B_i$
($1 \leq i \leq n)$. But $t \simeq\sb{\b} (x_r)v'_1\ldots v'_m$,
then $t \p_f (x_r)v_1\ldots v_m$ where $v_i \simeq\sb{\b} v'_i$
($1\leq i \leq m$), and, by induction hypothesis, \\$(u_r)v_1
[u_1/x_1,\ldots,u_n/x_n]\ldots v_{m-1} [u_1/x_1,\ldots,u_n/x_n]
\in |B\f
C|_{M,I}$ and \\$v_m [u_1/x_1,\ldots,u_n/x_n] \in |B|_{M,I}$. Therefore\\
$(u_r)v_1 [u_1/x_1,\ldots,u_n/x_n]\ldots v_m
[u_1/x_1,\ldots,u_n/x_n] \in |C|_{M,I}$ and\\ $t
[u_1/x_1,\ldots,u_n/x_n] \in |A|_{M,I}$.
\end{description}
\end{proof*}

\begin{definition}
We define the types with \textbf{positive quantifier} (resp.\
\textbf{ negative quantifier}) denoted $\q_2^+$ (resp.\ $\q_2^-$)
by:
 \begin{itemize}
\item An atomic formula is $\q_2^+$ and $\q_2^-$; \item If $A$ is
$\q_2^+$ (resp.\ $\q_2^-$) and $B$ is $\q_2^-$ (resp.\ $\q_2^+$),
then $B \rightarrow A$ is $\q_2^+$ (resp.\ $\q_2^-$); \item If $A$
is $\q_2^+$ and $x$ (resp.\ $X$) is a first order (resp.\ $n$-ary
second-order) variable, then $\q xA$ (resp.\ $\q XA$) is $\q_2^+$
; \item If $A$ is $\q_2^-$ and $x$ is a first-order variable, then
$\q xA$ is $\q_2^-$.
\end{itemize}
\end{definition}

We will prove that the $\q_2^+$  types are complete in ${{\cal
AF}2S}$.

\begin{definitions}
Let $\Omega = \{ x_i \,;\, i \in \textbf{N} \}$ be an enumeration
of an infinite set of variables of $\l$-calculus and $\{A_i \,;\,
i \in \textbf{N} \}$ be an enumeration of $\q_2^-$ types of
${{\cal AF}2S}$, where every $\q_2^-$ type occurs an infinite
number of times. We define the set $\mit \G^- = \{x_i : A_i \,; \,
i \in \textbf{N} \}$. Let $u$ be a $\l$-term such that $Fv (u)
\subseteq \Omega$, we define the contexte $\mit \G^-_u$ as the
restriction of $\mit \G^-$ on the set $Fv (u)$. The expression
$\mit \G^- \v_{{\cal AF}2S} u : B$ means that $\mit \G^-_u
\v_{{\cal AF}2S} u : B$. We put $\mit \G^- \v_{{\cal AF}2S}^{\b} u
: B$ iff there is a $\l$-term $u'$ such that $u \f_{\b} u'$ and
$\mit \G^- \v_{{\cal AF}2S} u' : B$.

Let $L$ be a second-order language and $E$ an equation system of
$L$. We define on the set of terms of $L$ an equivalence
 relation denoted $\approx_{E}$ by: $a \approx_{E} b$ iff we can
 obtain it by the following rules:
\begin{janclenum}{(\roman}{)}
\item if $a = b$ is a particular case of an equation of $E$, then
$a \approx_{E} b$;

\item for every terms $a, b, c$ of $L$, we have: $a \approx_{E}
 a$; and if $a \approx_{E} b$ and $b \approx_{E} c$, then $a
 \approx_{E} c$;

\item If $f$ is $n$-ary function symbol of $L$, and if $a_i
\approx_{E} b_i$ ($1 \leq i \leq n$), then $f (a_1,\ldots, a_n)
\approx_{E} f (b_1,\ldots, b_n)$.
\end{janclenum}
\end{definitions}

The following lemma allows to generalize the rule $(8)$.
\begin{Lem}\label{lemme:onze}
If $\G \v_{{\cal AF}2} u : B [a/x]$ and $a \approx_{E} b$, then
$\G \v_{{\cal AF}2} u : B [b/x]$.
\end{Lem}

\begin{proof*}
By induction in the definition of  $\approx_{E}$.
\end{proof*}

\begin{definition}
We consider $M_0$ the set of all closed terms of
$L$. We define a particular $\L_f$-model ${\cal M}$ by:
\begin{itemize}
\item The domain ${\cal |M|} = M_0 / \approx_{E}$ (the set of
equivalence classes modulo $\approx_{E}$); \item The adequate set
${\cal R}_f$;

\item To every $n$-ary symbol function $f$, we associate a
function
 $f_{\cal M} : |{\cal M}| ^n \f |{\cal M}|$ defined by
 $f_{\cal M} (\sou{a_1},\ldots, \sou{a_n}) = \sou{f
(a_1,\ldots, a_n)}$;

\item To every $n$-ary predicate symbol $P$, we associate a
function $P_{\cal M} : |{\cal M}| ^n \f {\cal R}_f$ defined by
$P_{\cal M} (\sou{a_1},\ldots, \sou{a_n}) = \{ \t \in \L \,;\,
\mit \G^- \v_{{\cal AF}2S}^{\b} \t : P(a_1,\ldots, a_n)\}$.
\end{itemize}
\end{definition}
It is easy to see that $f_{\cal M}$ and $P_{\cal M}$ are well defined.

\begin{definition}
We define a particular interpretation ${\cal I}$ on the variables
by: ${\cal I} (x) = \sou{x}$ and  ${\cal I} (X) = \Phi$, where
$\Phi : | {\cal M} |^n \f {\cal R}_f$ defined by $\Phi
(\sou{a_1},\ldots, \sou{a_n}) = \{ \t \in \L \,;\, \mit \G^-
\v_{{\cal AF}2S}^{\b} \t : X(a_1, \ldots, a_n)\}$.
\end{definition}

We have the following lemma.
\begin{Lem}\label{lemme:douze}
Let $S$ be a formula of $L$ and $\t$ a $\l$-term.
\begin{janclenum}{(\roman}{)}
\item If $S$ is $\q_2^+$ and $\t \in \mid S\mid_{\cal {M,I}}$,
then $\mit \G^-\v_{{\cal AF}2S}^{\b} \t : S$.\label{enum:i}

\item If $S$ is $\q_2^-$ and $\mit \G^-\v_{{\cal AF}2S}^{\b} \t :
S$, then $\t \in \mid S\mid_{\cal {M,I}}$.\label{enum:ii}
\end{janclenum}
\end{Lem}

\begin{proof*}
By simultanous induction on the $\q_2^+$ and $\q_2^-$ types.

\begin{subproof*}[of \textit{\ref{enum:i}\/}]{11.5cm}
\begin{enumerate}
\item $S$ is atomic: The result is trivial.

\item $S = \q XB$ where $B$ is $\q_2^+$: Let $Y$ be an $n$-ary
predicate variable which does not appear in $\mit \G_{\t}^-$ and
$B$. If $\t\in | \q XB |_{\cal {M,I}}$, then $\t\in | B[X]
|_{{{\cal M,I} [X \leftarrow | Y |_{\cal {M,I}}]}} = | B[Y / X]
|_{\cal {M,I}}$. By induction hypothesis, we have $\mit \G^-
\v_{{\cal AF}2S}^{\b} \t : B[Y]$, and there is a $\l$-term $\t'$
such that $\t \f_{\b} \t'$ and $\mit \G_{\t'}^- \v_{{\cal AF}2S}
\t' : B[Y]$. Since $Fv (\t') \subseteq Fv (\t)$, we deduce $\mit
\G_{\t'}^- \v_{{\cal AF}2S} \t' : \q YB[Y] = \q XB$, and $\mit
\G^- \v_{{\cal AF}2S}^{\b} \t : S$.

\item $S = B\f C$ where $B$ is $\q_2^-$ and $C$ is $\q_2^+$: let
$\t \in \mid B\f C \mid_{\cal {M,I}}$.  We put an $i$ such that $B
= A_i$ and $x_i$ is not free in $\t$. We have $x_i : B \v_{{\cal
AF}2} x_i : B$, then, by (ii), $x_i \in \mid B \mid_{\cal {M,I}}$,
therefore $(\t)x_i\in \mid C \mid_{\cal {M,I}}$, and, by induction
hypothesis, $\mit \G^- \v_{{\cal AF}2S}^{\b} (\t)x_i : C$. Thus
$(\t)x_i \f_{\b} \t'$ and $\mit \G^-_{\t'} \v_{{\cal AF}2S} \t' :
C$. We deduce that $(\t)x_i$ is normalizable, then $\t$ is also
normalizable. Since $(\t)x_i \f_{\b} \t'$, we obtain $\l
x_i(\t)x_i \f_{\b} \l x_i\t'$.

\begin{minipage}{11cm}
    \begin{itemize}
    \item
    If the normal form of $\t$ is $\l xu$, then $\l x_i(\t)x_i
    \f_{\b} \l x_i(\l xu)x_i \f_{\b} \l xu$ and $\l x_i\t'
    \f_{\b} \l xu$. But $\mit \G^- \v_{{\cal AF}2S} \l x_i\t' : S$ and $Fv
    (\l xu) \subseteq Fv (\l x_i\t')$, then, by Theorem
    \ref{theorem:cinq}, we obtain $\mit \G^- \v_{{\cal AF}2S} \l xu
    : S$, and $\mit \G^- \v_{{\cal AF}2S}^{\b} \t : S$.

    \item
    If not, let $v$ the normal form of $\t$. We have $\l
    x_i(\t)x_i \f_{\b} \l x_i(v)x_i$ and $\l x_i\t' \f_{\b} \l
    x_i(v)x_i$. Since $Fv (\l x_i(v)x_i) \subseteq Fv (\l x_i\t')$,
    we deduce that $\mit \G^- \v_{{\cal AF}2S} \l x_i(v)x_i : S$.
    Then, by Theorem \ref{theorem:six} and $Fv (\l x_i(v)x_i) = Fv (v)$,
    we obtain $\mit \G^- \v_{{\cal AF}2S} v : S$. Therefore $\mit \G^-
    \v_{{\cal AF}2S}^{\b} \t : S$.
    \end{itemize}
    \end{minipage}
\end{enumerate}
\end{subproof*}

\begin{subproof*}[of \textit{\ref{enum:ii}\/}]{11.5cm}
\begin{enumerate}
\item $S$ is atomic: The result is trivial. \item $S = B
\rightarrow C$ where $B$ is $\q_2^+$ and $C$ is $\q_2^-$: If $\mit
\G^- \v_{{\cal AF}2S}^{\b} \t : B \rightarrow C$, then, there is a
$\l$-term $\t'$ such that $\t \f_{\b} \t'$ and $\mit
\G^-_{\t'}\v_{{\cal AF}2S} \t' : B \rightarrow C$. If $u \in | B
|_{\cal {M,I}}$, then, by (i), $\mit \G^-\v_{{\cal AF}2S}^{\b}  u
: B$, and there is a $\l$-term $u'$ such that $u \f_{\b} u'$ and
$\mit \G^-_{u'}\v_{{\cal AF}2S} u' : B$. Therefore $\mit
\G^-_{(\t')u'}\v_{{\cal AF}2S} (\t') u' : C$, and, since $(\t)u
\f_{\b} (\t')u'$, we obtain $\mit \G^-\v_{{\cal AF}2S}^{\b} (\t)u
: C$. By induction hypothesis, we deduce $(\t)u \in | C |_{\cal
{M,I}}$.

\item $S = \q xB$ where $B$ is $\q_2^-$: Let $a \in | {\cal M} |$;
we
  have $a = \sou{b}$ where $b$ is a term of $L$. If $\mit \G^-
  \v_{{\cal AF}2S}^{\b} \t : \q xB$, then there is a $\l$-term $\t'$
  such that $\t \f_{\b} \t'$ and $\mit \G_{\t'}^-\v_{{\cal AF}2S} \t'
  : \q xB$, therefore $\mit \G_{\t'}^-\v_{{\cal AF}2S} \t' : B
  [b/x]$. But $B[b/x]$ is $\q_2^-$, then, by induction hypothesis,
  $\t' \in | B[b/x] |_{\cal {M,I}} = | B |_{{{\cal M,I} [x \leftarrow
    \sou{b}]}} = | B |_{{{\cal M,I} [x \leftarrow a]}}$. Thus $\t
  \in | B |_{{{\cal M,I} [x \leftarrow a]}}$ for every $a \in \cal
  |M|$.
\end{enumerate}
\end{subproof*}
\end{proof*}

\begin{The}\label{theorem:onze}
The closed $\q_2^+$ types are complete.
\end{The}

\begin{proof*}
Let $A$ be a closed $\q_2^+$ type. We will prove that: $t \in
|A|_{f}$ iff there is a $\l$-term $t'$ such that $t\f_{\b} t'$ and
$\v_{{\cal AF}2S} t' : A$.
 \begin{itemize}
\item That the condition is sufficient is a simple consequence of
Theorem \ref{theorem:dix}.

\item The condition is necessary: Indeed, let $t$ be a $\l$-term
such that $t \in | A |_{f}$, then $t \in |A|_{\cal M}$. We may
assume that $\mit \G_{t}^- = \emptyset$. By (i) of Lemma
\ref{lemme:douze}, we obtain $\mit \G^-\v_{{\cal AF}2S}^{\b} t :
A$, then there is a $\l$-term $t'$ such that $t \f_{\b} t'$ and
$\mit \G_{t'}^-\v_{{\cal AF}2S} t' : A$. Since $Fv (t') \subseteq
Fv (t)$, we deduce $\mit \G_{t'}^- = \emptyset$.
\end {itemize}
\end{proof*}

\begin{Cor}\label{corollaire:trois}
Let $A$ be a closed $\q_2^+$ type and $t$ a $\l$-term. If $t \in |A|_{f}$, then $t$
is normalizable and $\beta$-equivalent to a closed $\l$-term.
\end{Cor}

\begin{proof*}
By Theorem \ref{theorem:onze}.
\end{proof*}

\bibliography{biblio}

\end{document}